\documentclass[12pt,twoside]{article}

\title{A novel description and mathematical analysis of the Fractional Discrete Fourier Transform}

\author{Evan M. Zayas \\ Physics Department, Massachusetts Institute of Technology}
\date{\today}

\usepackage{amsmath}
\usepackage{amsfonts}
\usepackage{graphicx}
\usepackage{cancel}
\usepackage{indentfirst}
\usepackage[colorlinks=true]{hyperref}
\usepackage{array}
\usepackage{textcomp}
\usepackage[margin=1.5in]{geometry}


\begin{document}

\maketitle

\begin{abstract}
I discuss the nature of a Fractional Discrete Fourier Transform (FrDFT) described algorithmically by a combination of chirp transforms and ordinary DFTs. The transform is shown to be consistent with a continuous two-dimensional rotation between the time and frequency domains. I further present a new closed-form expression for the transformation matrix and some preliminary analysis of its properties.
\end{abstract}

\section{Motivation}

The Discrete Fourier Transform (DFT) has long been established as a principal tool for frequency analysis of signals. In particular, the transform decomposes a vector $\vec{x}$ in the time domain to its frequency spectrum $\vec{f}$ according to the transformation matrix:
\begin{align}
f_j &= \sum_{k=0}^{N-1} B_{jk} x_k \\
\label{bjk} \mathrm{with}\ &B_{jk} = \frac{1}{\sqrt{N}} e^{-\frac{2\pi i}{N}jk}\ \mathrm{for}\ 0 \leq j,k < N
\end{align}
 
\noindent where $N$ is the size of the data vector and $i \equiv \sqrt{-1}$. Furthermore, it is well understood that performing the (forward) DFT again on frequency-domain data returns it to the time domain but with the ordering reversed compared the original vector. We may verify this explicitly with the square of the transformation matrix above:
\begin{align}
\left(B^2\right)_{jk}=\sum_{l=0}^{N-1} B_{jl} B_{lk}&=\frac{1}{N}\sum_{l=0}^{N-1} e^{-\frac{2\pi i}{N}l(j+k)} \\
&=\begin{cases}1\ \mathrm{if}\ (j+k)\ \mathrm{mod}\ N=0 \\ 0\ \mathrm{else} \end{cases}
\end{align}

Thus, the matrix element is nonzero only when $k=N-j$ and $j>0$, or when instead $j=k=0$:
\begin{equation}
B^2=\begin{pmatrix} 1 & 0 & 0 & \dots \\ 0 & 0 & \dots & 1 \\ 0 & \dots & 1 & 0 \\ \vdots & \vdots & \vdots & \vdots \\ 0 & 1 & 0 & \dots \end{pmatrix}
\end{equation}

Acting with this matrix on an input vector $\vec{x}$ reverses the order of all elements except the first. It is common in such discussions to periodically extend the finite vector $\vec{x}$, so that $x_{j+N}=x_j$ for any integer $j$. Then, the action of two DFTs in series is equivalent to the parity flip operation $x_j \rightarrow x_{-j}$. Indeed, one can easily check from the above equation that $B^2$ is a square-root of the identity, consistent with the parity operator. This also immediately implies that the DFT itself is a fourth root of the identity.

From these properties, we can form a consistent picture of the Fourier Transform in the two-dimensional space of time and frequency. The DFT is interpreted as a rotation by $\pi/2$, for example from the time domain $(1, 0)$ to the orthogonal frequency domain $(0,1)$. Clearly, the subsequent interpretations of $B^2$ and $B^4$ which we have just discussed follow sensibly as well.

A monotone signal exists at a constant frequency for all time: a horizontal line $f(t)=f_0$ in this 2D description of the time and frequency domains. Its frequency spectrum is of course a delta function (in the large $N$ limit), which indicates that the signal is maximally localized in $f$; or equivalently, that a rotation of the coordinate space $(t,f)$ by $\pi/2$ localizes the signal in $t$. This description of the DFT provides a natural extension to the Fractional DFT (FrDFT), where a vector is instead rotated by an arbitrary angle $\alpha$. This makes accessible a continuum of mixed time/frequency domains parametrized by $\alpha$. The special property of a monotone signal is then extended to any line in the coordinate space: $f(t)=f_0+qt$, called a chirp signal. An appropriate choice of $\alpha$ fully localizes this signal just as the ordinary DFT does to a pure tone. The FrDFT therefore has potential applications in the study and processing of not only constant-frequency tones, but also chirped tones or locally-chirped tones including those reconstructed in experiments such as Project 8 \cite{p8} and LIGO \cite{ligo}.

In this work, I begin by studying an algorithm to perform the FrDFT proposed by Garcia {\it et al} \cite{garcia} in the context of this 2D rotation. I further use the Garcia algorithm to construct a simple closed-form expression for the FrDFT transformation matrix $F(\alpha)$. Lastly, I show that this expression for $F(\alpha)$ behaves consistently in the limiting cases of $\alpha=0$ and $\alpha=\pi/2$.

\section{Fractional DFT algorithm in time and frequency space}

The actual implementation of a FrDFT has been the subject of some considerable study for a variety of applications; of most interest to us is an algorithm published by Garcia {\it et al} \cite{garcia}. The algorithm consists only of two kinds of operations: (a) an ordinary DFT, and (b) multiplication by a quadratic phase. This is especially intriguing for two reasons: first, it is clearly related to chirp signals as the quadratic phase multiplication takes a plane wave $\exp{(-i \omega t)}$ to a chirp signal $\exp{(-i(\omega+qt) t)}$. Second, since the quadratic phase multiplication is an $\mathcal{O}(N)$ operation, this algorithm can be executed just as fast as the ordinary DFT, in $\mathcal{O}(N \log{N})$ time as is well established by Cooley and Tukey \cite{fft}. By contrast, direct computation of all matrix elements $F_{jk}(\alpha)$ is necessarily $\mathcal{O}(N^2)$ or worse.

The algorithm consists of five steps:

\begin{enumerate}
\item Multiply by the quadratic phase: $x_j \rightarrow \mathrm{exp}\left( -\frac{i\pi}{N} q_1 j^2 \right) x_j$
\item Perform the DFT
\item Multiply by the quadratic phase: $x_j \rightarrow \mathrm{exp}\left( -\frac{i\pi}{N} q_2 j^2 \right) x_j$
\item Perform the inverse DFT
\item Multiply again by the quadratic phase: $x_j \rightarrow \mathrm{exp}\left( -\frac{i\pi}{N} q_1 j^2 \right) x_j$
\end{enumerate}

\noindent where the chirp rates $q_1$ and $q_2$ are related to the rotation angle $\alpha$ by:
\begin{align}
q_1 &= \tan{(\alpha/2)}
\label{q1def} \\
q_2 &= \sin{(\alpha)}
\label{q2def}
\end{align}

To describe this algorithm in the $(t,f)$ coordinate space, we must find a 2x2 matrix that describes the action of each step; in fact, we have already done so for the DFT which takes $t \rightarrow f$ and $f \rightarrow -t$ via a $\pi/2$ rotation:
\begin{equation}
\left( \begin{array}{c} t' \\ f' \end{array} \right)=\left( \begin{array}{cc} 0 & -1 \\ 1 & 0 \end{array} \right) \left( \begin{array}{c} t \\ f \end{array} \right)
\end{equation}

The quadratic phase multiplication takes a constant tone with frequency $f$ to a chirp tone with frequency $f+qt$ and preserves $t \rightarrow t$:
\begin{align}
t'&=t \\
f'&=qt+f
\end{align}

\noindent and from the above equations we may simply read off the desired matrix elements. With the quadratic phase shift denoted $A(q)$ and the ordinary DFT denoted $B$, we have:
\begin{align}
A(q)&=\left( \begin{array}{cc} 1 & 0 \\ q & 1 \end{array} \right) \\
B&=\left( \begin{array}{cc} 0 & -1 \\ 1 & 0 \end{array} \right)
\end{align}

Note that with the above expression for $B$ we can once again easily check that the DFT is a fourth root of the identity. The FrDFT from the Garcia algorithm is:
\begin{equation}
\label{fa}
F(\alpha)=A(q_1)B^{-1}A(q_2)BA(q_1)
\end{equation}

\noindent and the claim is that with $q_1$ and $q_2$ related to $\alpha$ by Equations \ref{q1def} and \ref{q2def}, the matrix $F(\alpha)$ is a rotation by precisely the same angle $\alpha$. I will next explicitly show this claim to be true.
\begin{align}
\nonumber F(\alpha)&=\left( \begin{array}{cc} 1 & 0 \\ q_1 & 1 \end{array} \right)\left( \begin{array}{cc} 0 & 1 \\ -1 & 0 \end{array} \right)\left( \begin{array}{cc} 1 & 0 \\ q_2 & 1 \end{array} \right)\left( \begin{array}{cc} 0 & -1 \\ 1 & 0 \end{array} \right)\left( \begin{array}{cc} 1 & 0 \\ q_1 & 1 \end{array} \right) \\
&=\left( \begin{array}{cc} 1-q_1q_2 & -q_2 \\ q_1(2-q_1q_2) & 1-q_1q_2 \end{array} \right)
\end{align}

With Equation \ref{q2def} we can immediately replace $F_{01}(\alpha)=-q_2$ with $-\sin{(\alpha)}$. The diagonal elements both simplify to:
\begin{align}
\nonumber 1-q_1q_2&=1-\tan{(\alpha/2)}\sin{(\alpha)}\\
\nonumber &=1-\frac{\sin{(\alpha/2)}}{\cos{(\alpha/2)}}2\sin{(\alpha/2)}\cos{(\alpha/2)} \\
\nonumber &=1-2\sin^2{(\alpha/2)} \\
\Rightarrow 1-q_1q_2 &=\cos{(\alpha)}
\end{align}

\noindent and all that remains is to massage the last element $F_{10}$:
\begin{align}
\nonumber q_1(2-q_1q_2)&=q_2(2q_1/q_2-q_1^2) \\
\nonumber &=\sin{(\alpha)}\left(\frac{2 \tan{(\alpha/2)}}{2 \sin{(\alpha/2)}\cos{(\alpha/2)}}-\tan^2{(\alpha/2)}\right) \\
\nonumber &=\sin{(\alpha)}\left(\sec^2{(\alpha/2)}-\tan^2{(\alpha/2)}\right) \\
\Rightarrow q_1(2-q_1q_2)&=\sin{(\alpha)}
\end{align}

Indeed, we have shown that this representation of $F(\alpha)$ as a product of linear shifts and $\pi/2$ rotations is the general 2D rotation matrix:
\begin{equation}
F(\alpha)=\left( \begin{array}{cc} \cos{(\alpha)} & -\sin{(\alpha)} \\ \sin{(\alpha)} & \cos{(\alpha)} \end{array} \right)
\end{equation}

\section{Transformation matrix in the $N$-dimensional data space}

We have studied the Garcia algorithm in the coordinate space of time and frequency, and from it recovered the 2D rotation matrix. This solidifies our understanding of the DFT in this space, and provides good confidence that the FrDFT described by this algorithm is consistent with that understanding. Next, we return to the $N$-dimensional space of an input data vector to once again find an explicit expression for $F(\alpha)$.

The DFT in this space is already given by Equation \ref{bjk}, and we may similarly read off the matrix elements $A_{jk}(q)$ from the quadratic phase transformation:
\begin{align}
\tilde{x}_k=\mathrm{exp}\left(-\frac{i\pi}{N} q k^2\right) x_k&=\sum_{j=0}^{N-1} \mathrm{exp}\left(-\frac{i\pi}{N} q jk\right) x_j \delta_{jk} \\
\Rightarrow A_{jk}(q)&=\mathrm{exp}\left(-\frac{i\pi}{N}qjk\right) \delta_{jk}
\end{align}

\noindent where $\delta_{jk}$ is the Kronecker delta. For brevity, I now introduce the $2N^\mathrm{th}$ root of unity $\zeta=e^{-i\pi/N}$; the simplified expressions for $A_{jk}$ and $B_{jk}$ are:
\begin{align}
A_{jk}(q)&=\delta_{jk}\zeta^{qjk} \\
B_{jk}&=\frac{1}{\sqrt{N}}\zeta^{2jk}
\end{align}

\noindent and an explicit expression for $F_{jk}(\alpha)$ follows again from the matrix product as in Equation \ref{fa}:
\begin{align}
\nonumber F_{jk}&=\sum_{l,m,n,p=0}^{N-1} A_{jl}(q_1) B^{-1}_{lm} A_{mn}(q_2) B_{np} A_{pk}(q_1) \\
\nonumber &=\frac{1}{N}\sum_{l,m,n,p=0}^{N-1}\zeta^{q_1 jl} \zeta^{-2lm}\ \zeta^{q_2 mn}\zeta^{2np}\ \zeta^{q_1 pk}\ \delta_{jl} \delta_{mn} \delta_{pk} \\
\nonumber &=\frac{1}{N}\sum_{m=0}^{N-1} \zeta^{q_1 j^2} \zeta^{-2jm}\ \zeta^{q_2 m^2} \zeta^{2mk}\ \zeta^{q_1 k^2} \\
&=\frac{1}{N}\zeta^{q_1(j^2+k^2)}\sum_{m=0}^{N-1}\ \zeta^{q_2 m^2+2m(k-j)}
\end{align}

Substituting in $\zeta$ we have the full expression:
\begin{equation}
F_{jk}=\frac{1}{N}\ \mathrm{exp}-\frac{i\pi}{N}q_1(j^2+k^2)\sum_{m=0}^{N-1}\ \mathrm{exp}-\frac{i\pi}{N}\Big(q_2m^2+2m(k-j)\Big)
\label{fjk}
\end{equation}

This is perhaps a simpler closed-form expression for a FrDFT matrix than presently exists in the literature.

\section{Limiting behavior of $F(\alpha)$}
\label{lastSec}

Next, we will study this new expression for the FrDFT matrix in the two limiting cases of the rotation angle to prove the following results:

\begin{enumerate}
\item[(a)]{For $\alpha=0$, the transformation matrix $F_{jk}$ simply reduces to the identity: $F_{jk}(0)=\delta_{jk}$.}
\item[(b)]{For $\alpha=\pi/2$, the fractional transform corresponds to an ordinary DFT: $F_{jk}(\pi/2) = B_{jk}$ up to an overall phase.}
\end{enumerate}

It is quite straightforward to show that with $\alpha=0$ we recover the identity; both $q_1$ and $q_2$ are zero in this case, and thus $A_{jk}=\delta_{jk}$. The algorithm then reduces to a DFT followed by an inverse DFT, which trivially gives the identity. We can also arrive at this result directly from Equation \ref{fjk} with the $q_1$ and $q_2$ terms removed:
\begin{equation}
F_{jk}=\frac{1}{N} \sum_{m=0}^{N-1}\ \zeta^{2m(k-j)}=\delta_{jk}
\end{equation}

To prove (b), we set $\alpha=\pi/2$ which implies $q_1=q_2=1$ and simplify $F_{jk}$:
\begin{equation}
F_{jk}=\frac{1}{N} \zeta^{j^2+k^2} \sum_{m=0}^{N-1}\ \zeta^{m^2+2m(k-j)}
\end{equation}

Combining both factors of $\zeta$, we have in the exponent:
\begin{equation}
j^2+k^2+m^2+2m(k-j)=(m+k-j)^2+2jk
\end{equation}

The last term is precisely what we need to isolate $B_{jk}$:
\begin{equation}
F_{jk}=B_{jk}\frac{1}{\sqrt{N}} \sum_{m=0}^{N-1}\ \zeta^{(m+k-j)^2}
\label{fb}
\end{equation}

Thus, up to a normalization of $\sqrt{N}$ and an overall phase, we must only show that the sum in the above equation reduces to $1$, i.e. it is independent of $j$ and $k$. To accomplish this, I will first prove and then apply a trick with the roots of unity.

Consider a sum of the form:
\begin{equation}
S_k=\sum_{s=k}^{k+N-1} \zeta^{s^2}
\label{sdef}
\end{equation}

\noindent for arbitrary integer $k$. The difference between $S_{k+1}$ and $S_k$ is:
\begin{equation}
S_{k+1}-S_k=\zeta^{(k+N)^2}-\zeta^{k^2}
\end{equation}

\noindent since all of the remaining terms (when $k<s<N$) appear in both sums. We expand the exponent and recall that $\zeta^{2Ns}=1$ for any integer $s$, or equivalently that $\zeta^{Ns}=1$ for any even integer $s$. Then, provided only that $N$ is even, the difference between the sums is:
\begin{align}
\nonumber S_{k+1}-S_k&=\zeta^{k^2} \zeta^{2Nk} \zeta^{N^2}-\zeta^{k^2} \\
\nonumber &=\zeta^{k^2}-\zeta^{k^2} \\
\label{sk} \Rightarrow S_{k+1}=S_k
\end{align}

The result is a very general and somewhat surprising property of these sums: every $S_k$ is in fact the same. It follows by induction from Equation \ref{sk} that $S_k=S_0$ for every $k$, which explicitly removes the $k$ dependence:
\begin{equation}
\sum_{s=k}^{k+N-1} \zeta^{s^2}=\sum_{s=0}^{N-1} \zeta^{s^2}\ \text{for any integer }k
\label{trick}
\end{equation}

This completes the trick. Stated another way, we may freely displace the quantity $s$ in Equation \ref{sdef} by any integer $s \rightarrow s + \delta s$ and the sum $S_k=S_0$ remains unchanged. It is ultimately a simple consequence of the fact that $\zeta^{s^2}$ is periodic in $s$ with period $N$, and thus any sequence of $N$ terms will sum to the same value.

Returning now to Equation \ref{fb}, with this trick the majority of the proof becomes trivial. The remaining sum is exactly of the desired form with $s=m+k-j$, which we may simply replace with $s=m$:
\begin{equation}
F_{jk}=B_{jk}\frac{1}{\sqrt{N}} \sum_{m=0}^{N-1} \zeta^{m^2}
\end{equation}

\noindent and the $j,k$ dependence is now completely isolated to $B_{jk}$:
\begin{equation}
F_{jk}=\sigma\ B_{jk} \text{ with } \sigma=\frac{S_0}{\sqrt{N}}
\end{equation}

All that remains is to show $\sigma$ has unit norm. This is easily accomplished with the same trick, which we note is equally valid with $\zeta$ replaced by $\zeta^*$. The norm squared of $\sigma$ is:
\begin{equation}
|\sigma|^2=\frac{1}{N}\sum_{j=0}^{N-1} \zeta^{j^2}\ \sum_{k=0}^{N-1} \zeta^{-k^2}
\end{equation}

We displace the exponent in the sum over $k$ by $j$:
\begin{align}
\nonumber |\sigma|^2&=\frac{1}{N}\sum_{j=0}^{N-1} \zeta^{j^2}\ \sum_{k=0}^{N-1} \zeta^{-(k+j)^2} \\
\nonumber &=\sum_{k=0}^{N-1} \zeta^{-k^2}\left(\frac{1}{N} \sum_{j=0}^{N-1} \zeta^{-2jk}\right) \\
\nonumber &=\sum_{k=0}^{N-1} \zeta^{-k^2}\delta_{0k} \\
\Rightarrow |\sigma|^2&=1
\end{align}

Thus, up to an overall phase $\sigma=e^{i\theta}$ we have recovered the ordinary DFT matrix and obtained the desired result:
\begin{equation}
F_{jk}(\pi/2)=B_{jk}
\end{equation}

This completes the proof of (b).

\section{Concluding Remarks}

The work presented here has provided new insight into the nature of Fractional DFTs, but leaves a large amount open to potential future study. Using the Garcia algorithm, we have established a curious relationship between chirp transforms, rotations in two dimensions, and Fourier analysis. There is perhaps much more to be understood about the interconnected nature of these operations. Equation \ref{fjk} describes a FrDFT transformation matrix in closed form, from which each matrix element can be calculated from the sum of $N$ terms. This may also provide a new avenue to study the practical implementation of fractional transforms, and further analysis of this matrix like that in Section \ref{lastSec} to answer additional questions about its properties and other behavior.

\section*{Acknowledgements}

I would like to thank Sruthi Narayanan for insightful discussions pertaining to much of this work, and in particular for an initial proof of the roots of unity trick which was instrumental in reaching the conclusions of Section \ref{lastSec}. In addition, my research has been supported in part by NSF Award No. 1806251.


\begin{thebibliography}{99}

\bibitem{p8} A. A. Esfahani et al. {\it Determining the neutrino mass with cyclotron radiation emission spectroscopy$-$Project 8}. Journal of Physics G: Nuclear and Particle Physics {\bf 44}, 5. (2017)

\bibitem{ligo} E. Chassande-Motti et al. {\it Detection of GW bursts with chirplet-like template families}. Classical and Quantum Gravity {\bf 27}, 19. (2010)

\bibitem{garcia} J. Garcia et al. {\it Fractional-Fourier-transform calculation through the fast-Fourier-transform algorithm}. Applied Optics {\bf 35}, 35. (1996)

\bibitem{fft} J. W. Cooley and J. W. Tukey. {\it An algorithm for the machine calculation of complex Fourier series}. Mathematics of Computation {\bf 19}, pp. 297-301. (1965)

\end{thebibliography}
\end{document}